\documentclass[12pt]{article}
\usepackage{tikz}
\usetikzlibrary{calc,through,backgrounds}
\usepackage{multicol}

\usepackage{bm,amssymb,amsthm,amsfonts,amsmath,ytableau,latexsym,cite,color, psfrag,graphicx,ifpdf,pgfkeys,pgfopts,xcolor,extarrows}

\newtheorem{theo}{Theorem}[section]

\newtheorem{lem} [theo]{Lemma}

\newtheorem{prop}[theo]{Proposition}

\newtheorem{conj}[theo]{Conjecture}

\allowdisplaybreaks

\makeatletter \@addtoreset{equation}{section}
\@addtoreset{theo}{}\makeatother

\usepackage{hyperref}

\def\qed{\hfill \rule{4pt}{7pt}}
\def\pf{\noindent {\it Proof.} }
\def\S{  \mathfrak{S}}

\def\N{  \mathrm{N}}

\begin{document}
\begin{center}
{\bf COMPLEMENTS OF SCHUBERT POLYNOMIALS}

\vskip 4mm
{\small  NEIL J.Y. FAN, PETER L. GUO, NICOLAS Y. LIU}

\end{center}

\noindent{A{\scriptsize BSTRACT}.}
Let $\mathfrak{S}_w(x)$ be the Schubert polynomial for
a permutation  $w$ of $\{1,2,\ldots,n\}$. For any given composition $\mu$, we say that $x^\mu \mathfrak{S}_w(x^{-1})$ is
the complement of $\mathfrak{S}_w(x)$ with respect to $\mu$.
When each part of $\mu$ is  equal to $n-1$, Huh,   Matherne,   M\'esz\'aros and St.\,Dizier proved that  the normalization of $x^\mu \mathfrak{S}_w(x^{-1})$ is a Lorentzian polynomial.
They further conjectured that
the normalization of  $\mathfrak{S}_w(x)$ is Lorentzian. It can be shown  that if there exists a composition 
$\mu$ such that 
$x^\mu \mathfrak{S}_w(x^{-1})$ is a Schubert polynomial,
then  the normalization of  $\mathfrak{S}_w(x)$ will be Lorentzian. This motivates us to  investigate the problem of when  $x^\mu \mathfrak{S}_w(x^{-1})$
is  a Schubert polynomial.
We show that  if $x^\mu \S_w(x^{-1})$ is a Schubert polynomial, then $\mu$ must be a partition.
We also consider the   case when
 $\mu$ is the staircase partition $\delta_n=(n-1,\ldots, 1,0)$, and obtain  that  $x^{\delta_n} \S_w(x^{-1})$ is a Schubert polynomial if and only if $w$ avoids the   patterns 132 and 312. A conjectured characterization of when $x^\mu \S_w(x^{-1})$ is a Schubert polynomial is  proposed.


\vspace{6pt}


\section{Introduction}

The motivation of this paper is the Lorentzian property of Schubert polynomials  studied
 by Huh,   Matherne,   M\'esz\'aros and St.\,Dizier \cite{Huh}.
 For the beautiful theory developed around   Lorentzian  polynomials  as
 well as its powerful applications, see the
  work of  Br\"and\'en and   Huh \cite{Bra}.

Given  a polynomial $f=f(x_1,\ldots, x_n)\in \mathbb{C}[x_1,\ldots, x_n]$, the normalization of $f$,
denoted $\N(f)$, is the polynomial obtained from $f$ after applying  a linear operator $\N$ which is defined by
\[\N(x^\alpha)=\frac{x^\alpha}{\alpha_1!\cdots \alpha_n!},\ \ \ \ \ \text{for each monomial $x^\alpha=x_1^{\alpha_1}\cdots x_n^{\alpha_n}$}.\]
Let $S_n$ be the symmetric group of permutations of $[n]=\{1,2,\ldots,n\}$.
For $w\in S_n$, let $\S_w(x)=\S_w(x_1,\ldots,x_n)$ denote the Schubert polynomial
indexed by $w$. Schubert polynomials were introduced by Lascoux and Sch\"utzenberger \cite{Las}
as the polynomial representatives of Schubert classes in the  cohomology ring of the flag variety.

Using the Lorentzian property of volume polynomials of irreducible complex projective varieties, Huh,   Matherne,   M\'esz\'aros and St.\,Dizier \cite[Theorem 6]{Huh} showed that
\[
\N(x_1^{n-1}\cdots x_n^{n-1} \S_w(x^{-1})):=\N(x_1^{n-1}\cdots x_n^{n-1} \S_w(x_1^{-1},\ldots,x_n^{-1}))
\]
is  Lorentzian, see \cite{Bra,Huh} for the precise definition of Lorentzian polynomials.
As a consequence,  the normalization $\N(s_\lambda(x))$ of a Schur polynomial is Lorentzian \cite[Theorem 3]{Huh}. 
This  can be used to verify   Okounkov's log-concavity conjecture \cite[Conjecture 1]{Okounkov} for the Littlewood--Richardson coefficients  in
the   case of  Kostka numbers \cite[Theorem 2]{Huh}.

Huh,   Matherne,   M\'esz\'aros and St.\,Dizier \cite[Conjecture 15]{Huh} further
conjectured that  the normalization   $\N(\S_w(x))$ of each Schubert polynomial is
Lorentzian. Motivated by this conjecture, we   investigate the
following  problem:
\begin{align*}
\text{{\it When    $x^\mu \S_w(x^{-1})$
is still a Schubert polynomial?}}
\end{align*}
We call $x^\mu \S_w(x^{-1})$ the complement of $\S_w(x)$
with respect to a composition $\mu$. As will be seen in the following remark,
if there exists a  $\mu$ such that $x^\mu \S_w(x^{-1})$ is   a Schubert polynomial, then   $\N(\S_w(x))$ will satisfy the
Lorentzian property, and thus in this case we can verify the conjecture of
Huh,   Matherne,   M\'esz\'aros and St.\,Dizier.

\noindent
{\it Remark 1.} We explain that if there exists a $\mu$
such that $x^\mu \S_w(x^{-1})$ for $w\in S_n$ is   a Schubert polynomial, then   $\N(\S_w(x))$  is Lorentzian.
First, we claim  that   $\N(x^\mu\S_w(x^{-1}))$  is Lorentzian.
 This can be seen as follows.
It was observed  in \cite[Lemma 7]{Huh}   that for any $\nu\in \mathbb{Z}_{\geq 0}^n$ and a polynomial $f=f(x_1,\ldots, x_n)$, $\N(f)$ is Lorentzian if and only if $\N(x^\nu f)$ is Lorentzian.
Without loss of generality, we may assume that  $\max\{\mu_1,\ldots, \mu_n\}\leq n-1$ since
otherwise we  can embed $w$ into $S_m$ with $m>n$.
Let $\nu_i=n-1-\mu_i\in\mathbb{Z}_{\geq 0}$, so that
 \[x_1^{n-1}\cdots x_n^{n-1} \S_w(x^{-1})=x^\nu\cdot x^\mu\S_w(x^{-1}),\]
and hence,   the Lorentzian property  of $\N(x_1^{n-1}\cdots x_n^{n-1} \S_w(x^{-1}))$ \cite[Theorem 6]{Huh} implies that of  $\N(x^\mu\S_w(x^{-1}))$. This verifies the claim.
Suppose now that  $x^\mu \S_w(x^{-1})$  equals a Schubert polynomial, say $\S_{w'}(x)$.
Then we have $\S_w(x)=x^\mu \S_{w'}(x^{-1})$.
By the above claim,   $\N(x^\mu \S_{w'}(x^{-1}))$ is Lorentzian, and so it follows that
  $\N(\S_{w}(x))$  is Lorentzian.



Our first result shows that if $x^\mu \S_w(x^{-1})$ is a  Schubert polynomial,
then $\mu$ must be a partition, namely, the parts of $\mu$ are weakly decreasing.

\begin{theo}\label{mt1}
Let $w\in S_n$.
If there exists a composition  $\mu$ such that $x^\mu \S_w(x^{-1})$ is a Schubert polynomial, then $\mu$ is a partition.
\end{theo}

Let $\delta_n=(n-1,\ldots,1,0)$ denote the staircase partition. For $w\in S_n$, it is well known that if a monomial $x^\alpha$  appears in $\S_w(x)$, then  $x^\alpha|x^{\delta_n}$, that is, 
 $\alpha_i\leq n-i$ for $1\leq i\leq n$.
In  fact, the set $\{\S_w(x)\colon w\in S_n\}$ of Schubert polynomials
forms a $\mathbb{Z}$-basis of
\[V_n=\mathrm{span}_\mathbb{Z}\{x^\alpha\colon 0\leq \alpha_i\leq n-i\},\]
see for example \cite{Ber,Mac}.
Notice that
$x^\alpha\longmapsto x^{\delta_n-\alpha}$
induces    an involutive automorphism on $V_n$.
Hence the set  $\{x^{\delta_n} \S_w(x^{-1})\colon w\in S_n\}$   is a basis of $V_n$.
So it is natural  to investigate  the  typical case for $\mu=\delta_n$.

\begin{theo}\label{main1}
Let $w\in S_n$. Then  $x^{\delta_n} \S_w(x^{-1})$
is a Schubert polynomial if and only if
$w$ avoids the patterns $132$ and $312$.
\end{theo}

The sufficiency of Theorem \ref{main1} can be easily deduced  since
 the Schubert polynomial for a 132-avoiding permutation  is a single monomial.
However, the proof of the necessity of Theorem \ref{main1} turns out to
be quite technical.

\noindent
{\it Remark 2.}
The complements $x^{\delta_n} \S_w(x^{-1})$ are related to the  padded Schubert polynomials
defined by Gaetz and  Gao \cite{Gae,Gae-2}. The padded Schubert polynomial $\widetilde{\S}_w(x; y)$
is obtained from $\S_w(x)$ by replacing each monomial $x^\alpha$
 in $\S_w(x)$ by
$x^\alpha y^{\delta_n-\alpha}.$
Setting $x_i=1$ in $\widetilde{\S}_w(x; y)$, we see that $\widetilde{\S}_w(x; y)|_{x_i=1}=y^{\delta_n} \S_w(y^{-1})$. Thus  Theorem \ref{main1} is equivalent to saying that $\widetilde{\S}_w(x; y)|_{x_i=1}$ is a Schubert
polynomial if and only if $w$ avoids the patterns $132$ and $312$.

Computational  verification suggests the following conjecture.

\begin{conj}\label{cytj}
Let $w\in S_n$.
Then there exists  a partition $\mu$ such that $x^\mu \S_w(x^{-1})$ is a Schubert polynomial
if and only if $w$ avoids the  patterns  1432, 13254, 14253, 24153, 31524, 361452.
\end{conj}

This paper is organized as follows. In Section \ref{ma1}, we give a proof of Theorem \ref{mt1}.
The arguments rely on the RC-graph (also called pipedream) model  of Schubert polynomials.
In Section \ref{RCC}, we present a proof of Theorem \ref{main1}.
Remarks about Conjecture \ref{cytj} are given in Section \ref{TEE}.

\vspace{.2cm} \noindent{\bf Acknowledgments.}
We would like to thank the anonymous referee for valuable comments and suggestions. 
This work was
supported by the National Science Foundation of China (11971250, 12071320, 12371329)
and Sichuan Science and Technology Program (2023ZYD0012).

\section{Proof of Theorem \ref{mt1}}\label{ma1}

To give a proof of Theorem \ref{mt1}, we require the RC-graph model of Schubert polynomials
investigated by Bergeron and Billey \cite{Ber}.

Schubert polynomials are defined based on divided difference operators .
For a polynomial $f(x)\in \mathbb{Z}[x_1,\ldots,x_n]$, the  divided difference operator $\partial_i$
acts on  $f(x)$ by  letting
\[\partial_i f(x)=\frac{f(x)
    -s_if(x)}{x_i-x_{i+1}},\]
where $s_if(x)$ is obtained from $f(x)$ by exchanging $x_i$ and $x_{i+1}$.
It is easily checked that  $\partial_if(x)$  is symmetric in $x_i$ and $x_{i+1}$.

For $w_0=n \cdots   21$, set
   \[\S_{w_0}(x)=x_1^{n-1}x_2^{n-2}\cdots x_{n-1}.\]
 For $w=w_1w_2\cdots w_n\neq w_0$, choose a position $i$ such that $w_i<w_{i+1}$ and define
\[
\S_{w}(x)=\partial_i \S_{ws_i}(x),
\]
where $ws_i$ is   obtained from $w$ by  interchanging   $w_i$ and $w_{i+1}$.
The  above definition is independent of the choice of $i$ since the operators $\partial_i$ satisfy the
Coxeter relations: $\partial_i \partial_j=\partial_j \partial_i$ for $|i-j|>1$, and
$\partial_i \partial_{i+1} \partial_i= \partial_{i+1}\partial_i  \partial_{i+1}$.

 Billey, Jockusch and  Stanley \cite{Bil} showed that  $\S_{w}(x)$ can be combinatorially generated 
  in terms of reduced-word compatible sequences.
Each permutation  $w$ can be expressed as a product  $s_{a_1}s_{a_2}\cdots s_{a_p}$ of
  adjacent transpositions. If $p$ is minimal, then we write   $\ell(w)=p$,
  which is   the length of $w$, and in this case
 $a=(a_1,\ldots, a_p)$  is called a reduced word of $w$.
 A sequence $\alpha=(\alpha_1,\ldots, \alpha_p)$
 of positive integers is said to be a compatible sequence of $a$ if (i)
 $\alpha_1\leq \alpha_2\leq \cdots\leq \alpha_p$, (ii) $\alpha_i< \alpha_{i+1}$ if $a_i<a_{i+1}$,
 and (iii) $1\leq \alpha_i\leq a_i$.
 Let $\mathrm{Red}(w)$ denote the set of reduced words of $w$.
 For a reduced word $a\in \mathrm{Red}(w)$, $\mathrm{C}(a)$ denotes the set of its compatible sequences.
It was proved in  \cite{Bil}  that
for $w\in S_n$,
\begin{equation}\label{BJS}
\S_w(x)=\sum_{a\in \mathrm{Red}(w)}\sum_{\alpha\in \mathrm{C}(a)}x_{\alpha_1}x_{\alpha_2}\cdots x_{\alpha_p}.
\end{equation}

The RC-graph
corresponding to a compatible pair $(a,\alpha)$ can be defined as follows.
Let  $\Delta_n$ denote the array of
 left-justified boxes  with $n+1-i$ boxes in row $i$ for $1\leq i\leq n$.
 The RC-graph
corresponding to  $(a,\alpha)$ is the subset  of $\Delta_n$ consisting of boxes in row $\alpha_i$
 and column $a_i-\alpha_i+1$, where $1\leq i\leq p$. For example, let  $a=(4,2,3,2,4)$ be a
 reduced word of $w=15342$, and $\alpha=(1,1,2,2,4)$ be a compatible sequence of $a$.
 Then the corresponding RC-graph   is illustrated  in Figure \ref{rc-graph},
 where we use a cross to signify a box belonging to the RC-graph.
 \begin{figure}[h]
\setlength{\unitlength}{0.5mm}
\begin{center}
\begin{picture}(50,50)
\put(0,-5){\line(1,0){10}}\put(0,5){\line(1,0){20}}\put(0,15){\line(1,0){30}}
\put(0,25){\line(1,0){40}}\put(0,35){\line(1,0){50}}\put(0,45){\line(1,0){50}}

\put(50,45){\line(0,-1){10}}\put(40,45){\line(0,-1){20}}\put(30,45){\line(0,-1){30}}
\put(20,45){\line(0,-1){40}}\put(10,45){\line(0,-1){50}}\put(0,45){\line(0,-1){50}}

\put(2,8){$+$}\put(2,28){$+$}\put(12,28){$+$}\put(12,38){$+$}\put(32,38){$+$}

\end{picture}
\end{center}\vspace{-.2cm}
\caption{An RC-graph of $w=15342$.}
\label{rc-graph}
\end{figure}

For an RC-graph $D$ of $w$, let
\[x^D=\prod_{i\in [n]} x_i^{|\text{\{boxes in the $i$-th row of $D$\}}|}.\]
In this notation, formula \eqref{BJS} can be rewritten as
\begin{equation*}
\S_w(x)= \sum_{D} x^D,
\end{equation*}
where the sum runs over all the RC-graphs of $w$.

For the purpose of this paper, we pay attention to two specific RC-graphs:
the bottom RC-graph and the top RC-graph.
Let $d(w)=(d_1,d_2,\ldots, d_n)$ denote the inversion code of $w$, that is,
for $1\leq i\leq n$,
 \[d_i=|\{j\colon j>i,\ w_j<w_i\}|.\]
 Clearly, $0\leq d_i\leq n-i$.
  The bottom RC-graph of $w$  is the RC-graph  consisting of the first $d_i$ boxes in row $i$.
For example, Figure \ref{BT}(a) is the bottom RC-graph of $w=25143$.
 \begin{figure}[h]
\setlength{\unitlength}{0.5mm}
\begin{center}
\begin{picture}(150,50)
\put(0,0){\line(1,0){10}}\put(0,10){\line(1,0){20}}\put(0,20){\line(1,0){30}}
\put(0,30){\line(1,0){40}}\put(0,40){\line(1,0){50}}\put(0,50){\line(1,0){50}}

\put(50,50){\line(0,-1){10}}\put(40,50){\line(0,-1){20}}\put(30,50){\line(0,-1){30}}
\put(20,50){\line(0,-1){40}}\put(10,50){\line(0,-1){50}}\put(0,50){\line(0,-1){50}}

\put(2,13){$+$}\put(2,33){$+$}\put(12,33){$+$}\put(22,33){$+$}
\put(2,43){$+$}


\put(100,0){\line(1,0){10}}\put(100,10){\line(1,0){20}}\put(100,20){\line(1,0){30}}
\put(100,30){\line(1,0){40}}\put(100,40){\line(1,0){50}}\put(100,50){\line(1,0){50}}

\put(150,50){\line(0,-1){10}}\put(140,50){\line(0,-1){20}}\put(130,50){\line(0,-1){30}}
\put(120,50){\line(0,-1){40}}\put(110,50){\line(0,-1){50}}\put(100,50){\line(0,-1){50}}

\put(102,33){$+$}\put(122,33){$+$}\put(122,43){$+$}
\put(102,43){$+$}\put(132,43){$+$}

\put(15,-10){(a)}\put(115,-10){(b)}
\end{picture}
\end{center}
\caption{The bottom and   top RC-graphs of $w=25143$.}
\label{BT}
\end{figure}

Bergeron and Billey \cite{Ber}  showed that
any RC-graph of $w$ can be obtained  from the bottom RC-graph of $w$ by
applying a sequence of ladder moves. Let $D$ be an RC-graph   of
$w$, and let $(i,j)$ denote a box of $D$ in row $i$ and $j$. The ladder move $L_{i,j}$  is a
local change of the crosses  as illustrated in Figure \ref{ladder}.
Formally, the resulting diagram after applying the ladder move $L_{i,j}$ is
\[L_{i,j}(D)=D\backslash\{(i,j)\}\cup\{(h,j+1)\}.\]
 \begin{figure}[h]
\setlength{\unitlength}{0.5mm}
\begin{center}
\begin{picture}(100,60)
\put(0,-5){\line(1,0){20}}\put(0,5){\line(1,0){20}}\put(0,15){\line(1,0){20}}
\put(0,30){\line(1,0){20}}\put(0,40){\line(1,0){20}}\put(0,50){\line(1,0){20}}

\put(0,-5){\line(0,1){55}}\put(10,-5){\line(0,1){55}}\put(20,-5){\line(0,1){55}}

\put(2,-2){$+$}\put(2,8){$+$}\put(12,8){$+$}\put(2,33){$+$}
\put(12,33){$+$}\put(4,19){$\vdots$}\put(14,19){$\vdots$}


\put(80,-5){\line(1,0){20}}\put(80,5){\line(1,0){20}}\put(80,15){\line(1,0){20}}
\put(80,30){\line(1,0){20}}\put(80,40){\line(1,0){20}}\put(80,50){\line(1,0){20}}

\put(80,-5){\line(0,1){55}}\put(90,-5){\line(0,1){55}}\put(100,-5){\line(0,1){55}}

\put(92,43){$+$}\put(82,8){$+$}\put(92,8){$+$}\put(82,33){$+$}
\put(92,33){$+$}\put(84,19){$\vdots$}\put(94,19){$\vdots$}

\put(35,28){$\underrightarrow{\quad L_{i,j}\quad }$}

\put(-8,-2){{\small $i$}}\put(3,54){{\small $j$}}
\put(72,-2){{\small $i$}}\put(83,54){{\small $j$}}
\put(-8,43){{\small $h$}}\put(72,43){{\small $h$}}
\end{picture}
\end{center}
\caption{A ladder move.}
\label{ladder}
\end{figure}
It can be  verified that $L_{i,j}(D)$ is still an
RC-graph of $w$.
 Since the ladder move operation
moves   a cross  upwards, the bottom RC-graph of $w$
corresponds to the leading monomial of $\S_w(x)$ in the reverse lexicographic order.

\begin{prop}[\mdseries{Bergeron--Billey \cite{Ber}}]\label{large}
Let $w\in S_n$, and $d(w)=(d_1,\ldots,d_n)$ be the inversion code of $w$.
Then the monomial  $x^{d(w)}=x_1^{d_1}\cdots x_n^{d_n}$ is the leading term
of $\S_w(x)$ in the reverse lexicographic order.
\end{prop}

It was shown in \cite{Ber} that  the transpose of  an RC-graph  of $w$
is an RC-graph of the inverse  $w^{-1}$ of $w$.  The top RC-graph of $w$  is
defined as the transpose of the bottom RC-graph of $w^{-1}$. For example, the inverse of
$w=25143$ is $w^{-1}=31542$, and so the top RC-graph of $w$ is as depicted in Figure \ref{BT}(b).

Dual to the ladder moves, Bergeron and Billey \cite{Ber} defined  chute moves  on
RC-graphs, and showed that
any RC-graph of $w$ can be obtained  from the top RC-graph of $w$ by
applying a sequence of chute moves.
The chute move operation  moves a cross downwards, and
thus the top RC-graph  of $w$
corresponds to the smallest monomial of $\S_w(x)$ in the reverse lexicographic order.

For a vector   $v=(v_1,v_2,\ldots,v_n)$ with $0\leq v_i\leq n-i$,
let
\begin{equation}\label{trans}
v^t=(v_1', v_2',\ldots, v_n')
\end{equation}
 denote the transpose of $v$,
namely, for $1\leq i\leq n$,
\[v_i'=|\{1\leq j\leq n\colon v_j\geq i\}|.\]
It is easy to see that $v^t$ is a partition for any $v$.

\begin{prop}[\mdseries{Bergeron--Billey \cite{Ber}}]\label{small}
For $w\in S_n$, let
$d^{\,t}(w^{-1})=(d_1',\ldots, d_n')$
denote the transpose of the inversion code
   of $w^{-1}$. Then the monomial
$x^{d^{\,t}(w^{-1})}=x_1^{d_1'}\cdots x_n^{d_n'}$
 is the smallest term
of $\S_w(x)$ in the reverse lexicographic order.
\end{prop}

Using  Propositions \ref{large} and \ref{small}, we
can now provide a proof of  Theorem \ref{mt1}.

\noindent{\it Proof of Theorem \ref{mt1}.}
 Suppose that
\begin{equation}\label{QQQ1}
x^\mu \S_w(x^{-1})=\S_{w'}(x).
\end{equation}
Equivalently,
\begin{equation}\label{QQQ2}
x^\mu \S_{w'}(x^{-1})=\S_{w}(x).
\end{equation}
  By Propositions \ref{large} and \ref{small}, it follows that
\[d(w')=\mu-d^t(w^{-1})\ \ \ \ \ \text{and}\ \ \ \ \ d(w)=\mu-d^t(w'^{-1}),\]
 and so we have
\begin{equation}\label{NBN}
d(w')+d^t(w^{-1})=d(w)+d^t(w'^{-1})=\mu.
\end{equation}
 Suppose to the contrary that $\mu$ is not a partition. Locate  an index $i$ such that $\mu_i<\mu_{i+1}$. Write $d(w)=(d_1,d_2,\ldots, d_n)$ and $d(w')=(d_1', d_2',\ldots, d_n')$. Keep in mind that both $d^t(w^{-1})$ and $d^t(w'^{-1})$ are partitions. In view of \eqref{NBN} together with the assumption $\mu_i<\mu_{i+1}$,
  we must have $d_i<d_{i+1}$ and $d_i'<d_{i+1}'$. By the definition of inversion code, we obtain that $w_i<w_{i+1}$ and $w'_i<w'_{i+1}$. Therefore,
  \[\S_w(x)=\partial_i \S_{ws_i}\ \ \ \ \ \text{and}\ \ \ \ \ \S_{w'}(x)=\partial_i \S_{w's_i},\]
 implying that both $\S_w(x)$ and $\S_{w'}(x)$ are symmetric in $x_i$ and $x_{i+1}$.
Since $\mu_i<\mu_{i+1}$,   both $x^\mu \S_w(x^{-1})$ and $x^\mu \S_{w'}(x^{-1})$ could not be symmetric in $x_i$ and $x_{i+1}$, which clearly contradicts \eqref{QQQ1} and \eqref{QQQ2}. This completes the proof.
\qed

\section{Proof of Theorem \ref{main1}}\label{RCC}

As mentioned in Introduction, the main difficulty  is to prove the necessity of
Theorem \ref{main1}, that is,
if $x^{\delta_n} \S_w(x^{-1})$
is a Schubert polynomial, then
$w$ must avoid  the patterns $132$ and $312$.

\subsection{Permutations avoiding 132 and 312 }

In this subsection, we prove a relationship  concerning   permutations avoiding the
patterns $132$ and $312$, see Proposition \ref{lemma-a}. This proposition
will be used in the   proof of Lemma \ref{LL-4}, which is crucial to  prove the necessity of
Theorem \ref{main1}.

Given a permutation $w=w_1w_2\cdots w_n\in S_n$, we say that $w$ is 132-avoiding
 if there do not exist $i<j<k$ such that $w_i<w_k<w_j$, and $w$ is 312-avoiding
if there do not exist $i<j<k$ such that $w_j<w_k<w_i$.
For $1\leq i\leq n$, define
\begin{align}\label{lr}
l_i(w)=|\{j\colon j<i,w_j<w_{i}\}| \ \ \ \text{and}\ \ \
r_i(w)=|\{j\colon j<i,w_j>w_{i}\}|.
\end{align}
Clearly, $l_i(w)+r_i(w)=i-1$.
Let
\begin{align}
L(w)=(l_1(w),\ldots,l_n(w)) \ \ \ \text{and}\ \ \ R(w)=(r_1(w),\ldots,r_n(w)).
\end{align}

We say that a vector  $v=(v_1,\ldots,v_n)$ is a rearrangement of
a vector $v'=(v'_1,\ldots,v'_n)$  if there exists a permutation $\pi\in S_n$
such that $v=(v'_{\pi_1},\ldots,v'_{\pi_n})$.

\begin{prop}\label{lemma-a}
Let $w\in S_n$ be a  permutation avoiding 312, and $u\in S_n$
be a permutation avoiding 132. If $L(w)$ is a rearrangement of
$L(u)$ and $R(w)$ is a rearrangement of
$R(u)$, then we have $w=u$.
\end{prop}

To prove Proposition \ref{lemma-a}, we need the matrix representation  of $w$.
Consider  an $n\times n$ box grid, where the rows (respectively, the columns) are numbered
 $1,2,\ldots,n$ from
top to bottom (respectively, from left to right). A box in row $i$
and column $j$ is denoted  $(i,j)$. The  matrix representation  of $w$ is obtained
by putting a dot in the box $(i,w_i)$ for $1\leq i\leq n$. Figure \ref{Rothe}
is the matrix representation  of $w=426315$.
\begin{figure}[h]
\setlength{\unitlength}{0.5mm}
\begin{center}
\begin{picture}(60,65)
\put(0,0){\line(1,0){60}}\put(0,10){\line(1,0){60}}\put(0,20){\line(1,0){60}}
\put(0,30){\line(1,0){60}}\put(0,40){\line(1,0){60}}\put(0,50){\line(1,0){60}}
\put(0,60){\line(1,0){60}}

\put(0,0){\line(0,1){60}}\put(10,0){\line(0,1){60}}
\put(20,0){\line(0,1){60}}
\put(30,0){\line(0,1){60}}\put(40,0){\line(0,1){60}}
\put(50,0){\line(0,1){60}}
\put(60,0){\line(0,1){60}}

\put(55,35){\circle*{3}}\put(45,5){\circle*{3}}
\put(25,25){\circle*{3}}\put(5,15){\circle*{3}}
\put(35,55){\circle*{3}}
\put(15,45){\circle*{3}}

\end{picture}
\end{center}
\vspace{-.5cm}
\caption{ The  matrix representation  of $w=426315$.}
\label{Rothe}
\end{figure}
 Obviously,  for $1\leq i\leq n$, the value $l_i(w)$ (respectively, $r_i(w)$) is equal to the
number of dots lying in the region to the strictly upper
left (respectively, strictly  upper right) of the box $(i, w_i)$.
For $w=426315$, it can be seen from Figure \ref{Rothe} that
$L(w)=(0,0,2,1,0,4)$ and $R(w)=(0,1,0,2,4,1).$

\vspace{10pt}

\noindent
{\it Proof of Proposition \ref{lemma-a}.}
We first prove that either $w_n=u_n=n$ or $w_n=u_n=1$. Since $w$ is 312-avoiding, we see that
for $1\leq j,\, k\leq n-1$, if $w_j<w_n$ and $w_k>w_n$, then we must have $j<k$. This implies
that the matrix representation of
$w$ is as illustrated in Figure \ref{pf-2}, where the region $L$ contains all the $l_n(w)$ dots
that are on the upper left of $(n, w_n)$, and the
region $R$ contains all the $r_n(w)=n-1-l_n(w)$ dots
that are on the upper right of $(n, w_n)$.
We have the following observation:
\begin{itemize}
\item[(O1).]  For $1\leq j \leq l_n(w)$, we have $l_j(w)< l_n(w)$, and for $l_n(w)< j\leq n$,
we have $l_j(w)\geq l_n(w)$. Therefore, there are exactly
\[n-l_n(w)=r_n(w)+1\]
 indices
$j$ such that $l_j(w)\geq l_n(w)$.
\end{itemize}

\begin{figure}[h]
\setlength{\unitlength}{0.5mm}
\begin{center}
\begin{picture}(40,50)
\put(0,0){\line(1,0){50}} \put(0,10){\line(1,0){50}}
\put(0,30){\line(1,0){50}}\put(0,50){\line(1,0){50}}

\put(0,0){\line(0,1){50}} \put(20,0){\line(0,1){50}}
\put(30,0){\line(0,1){50}} \put(50,0){\line(0,1){50}}

\put(-15,50){\line(1,0){13}}\put(-15,30){\line(1,0){13}}
\put(52,30){\line(1,0){13}}
\put(52,10){\line(1,0){13}}

\put(-8,50){\vector(0,-1){4}}\put(-8,30){\vector(0,1){4}}
\put(-15,39){\tiny{$l_n(w)$}}

\put(58,30){\vector(0,-1){4}}\put(58,10){\vector(0,1){4}}
\put(52,19){\tiny{$r_n(w)$}}

\put(25,5){\circle*{3}}

\put(7,37){$L$}\put(37,17){$R$}

\put(-8,3.5){\small{$n$}}\put(20.5,53){\small{$w_n$}}

\end{picture}
\end{center}
\vspace{-.5cm}
\caption{ The  matrix representation  of $w$.}
\label{pf-2}
\end{figure}

Since $l_n(w)+r_n(w)=n-1$,
we have the following two cases: $l_n(w)\geq \frac{n-1}{2}$ and $r_n(w)>\frac{n-1}{2}$.

Case 1. $l_n(w)\geq \frac{n-1}{2}$. Since $L(w)$ is a rearrangement of $L(u)$, there
exists some $1\leq i\leq n$ such that $l_i(u)=l_n(w)$, that is, in the matrix
representation of $u$,
 there are $l_n(w)$
dots lying in the region to the upper left of the box $(i, u_i)$. Since $u$ is 132-avoiding,
we see that for $1\leq j, k<i$, if $u_j<u_i$ and $u_k>u_i$,
 then  $j>k$. So the matrix
representation  of $u$ is as illustrated in Figure \ref{pf-1},
where the box  marked with a dot is
$(i,u_i)$, the region $L'$ contains  $l_i(u)=l_n(w)$ dots, and
the region $A$ contains  $r_i(u)$ dots.
\begin{figure}[h]
\setlength{\unitlength}{0.5mm}
\begin{center}
\begin{picture}(50,50)
\put(0,0){\line(1,0){50}} \put(0,20){\line(1,0){50}}
\put(0,30){\line(1,0){50}}\put(0,40){\line(1,0){50}}\put(0,50){\line(1,0){50}}

\put(0,0){\line(0,1){50}} \put(20,0){\line(0,1){50}}
\put(30,0){\line(0,1){50}} \put(50,0){\line(0,1){50}}

\put(25,25){\circle*{3}}\put(-5,22){\small{$i$}}
\put(21.5,53){\small{$u_i$}}

\put(7,32){$L'$}\put(37,42){$A$}
\put(7,7){$B$}\put(37,7){$C$}
\end{picture}
\end{center}
\vspace{-.5cm}
\caption{ The  matrix representation  of $u$ in Case 1.}
\label{pf-1}
\end{figure}

We next show that there are no dots in the regions $A, B, C$ of $u$ in Figure \ref{pf-1}.
To proceed, we list the following straightforward  observations.
\begin{itemize}
\item[(O2).] Since there are $l_n(w)\geq \frac{n-1}{2}$ dots in   $L'$,
 there are at most  $\frac{n-1}{2}$ dots in  $A$. So,
 if $(j, u_j)$ is a box in  $A$, then $l_j(u)<l_n(w)$.

\item[(O3).]  Since $u$ is 132-avoiding,  every dot (if any) in   $B$
lies to the left of  every dot in  $L'$. Again,
since there are $l_n(w)\geq \frac{n-1}{2}$ dots in  $L'$, there are
at most  $\frac{n-1}{2}$ dots in  $B$. So
 if $(j, u_j)$ is a box in  $B$, then $l_j(u)<l_n(w)$.

\item[(O4).]  If $(j, u_j)$ is a box in   $C$, then $l_j(u)\geq l_n(w)$.

\end{itemize}

By  (O2),  (O3) and (O4), we see that the set of indices $j$
such that $l_j(u)\geq l_n(w)$ is
\[\{i\} \cup \{k\colon \text{$(k, u_k)$ is a box in   $C$}\}.\]
Keep in mind that
$L(w)$ is a rearrangement of $L(u)$. By    (O1), we  have
\begin{equation*}
|\{i\} \cup \{k\colon \text{$(k, u_k)$ is a box in  $C$}\}|=1+r_n(w),
\end{equation*}
and so
\begin{equation}\label{XX-1}
|\{\text{dots in   $C$}\}|=r_n(w).
\end{equation}
Notice that the total number of dots in  $A$, $B$ and $C$
is equal to $n-1-l_n(w)=r_n(w)$.
This together with \eqref{XX-1} forces that there are no dots
in  $A$ and $B$.

It remains to show that there are  no dots in  $C$.
Suppose otherwise there is at least one dot in $C$.
Note that the dot in $(i, u_i)$ appears to the
upper left of each dot in  $C$. Since $u$ is 132-avoiding,  the dots in  $C$ must be
 listed
from upper left to bottom right increasingly.
This, along with the fact that there are no dots in  $A$ and $B$, implies
that $u_n=n$. So we have $l_n(u)=n-1$. Since $L(w)$ is a rearrangement of
$L(u)$, there exists some index, say  $k$, such that $l_k(w)=n-1$.
By the definition of $l_k(w)$, we must have $k=n$ and  $w_n=n$.
Recalling that $l_i(u)=l_n(w)$,  we have $l_i(u)=l_n(w)=n-1$,
and thus $i=n$.
However,  from the assumption that
 there is at least one dot in $C$, it follows that $i<n$,
leading to a contradiction. Hence the region $C$ does not contain any dots.

Because there are no dots in $B$ and $C$, we have $i=n$. By \eqref{XX-1}, we see that $r_n(w)=0$,
and so we have $w_n=n$.
 Moreover, since $l_n(u)=l_n(w)=n-1-r_n(w)=n-1$,
we see that $u_n=n$. Therefore, we have $w_n=u_n=n$.

Case 2. $r_n(w)>\frac{n-1}{2}$. In this case, from   Figure \ref{pf-2},
we see that $r_n(w)$ is the unique maximum value in $R(w)$.
Since $R(u)$ is a rearrangement of $R(w)$, there is a unique $i$ such that
$r_i(u)=r_n(w)$. The matrix representation of $u$ is illustrated in
Figure \ref{pf-3}, where the square marked with a dot is $(i, u_i)$ and the region
$R'$ contains $r_i(u)=r_n(w)$ dots.
\begin{figure}[h]
\setlength{\unitlength}{0.5mm}
\begin{center}
\begin{picture}(50,50)
\put(0,0){\line(1,0){50}} \put(0,20){\line(1,0){50}}
\put(0,30){\line(1,0){50}}\put(0,40){\line(1,0){50}}\put(0,50){\line(1,0){50}}

\put(0,0){\line(0,1){50}} \put(20,0){\line(0,1){50}}
\put(30,0){\line(0,1){50}} \put(50,0){\line(0,1){50}}

\put(25,25){\circle*{3}}\put(-5,22){$i$}\put(21,53){$u_i$}

\put(7,32){$D$}\put(37,42){$R'$}
\put(7,7){$E$}\put(37,7){$F$}
\end{picture}
\end{center}
\vspace{-.5cm}
\caption{ The  matrix representation  of $u$ in Case 2.}
\label{pf-3}
\end{figure}

Since $r_i(u)=r_n(w)>\frac{n-1}{2}$ is the unique maximum value in $R(u)$, we see that
there are no dots in the regions $D$ and $E$. We next show that
the region $F$   does not contain any dots either.
By the fact that $l_n(w)+r_n(w)=n-1$,
as well as the fact that  $D$ and $E$ do not contain dots,
it follows  that there are $l_n(w)$ dots in $F$:
\begin{equation}\label{buch}
 |\{\text{dots in $F$}\}|=l_n(w).
 \end{equation}
We aim to verify $l_n(w)=0$.
Suppose otherwise that $l_n(w)>0$. By   (O1),
there are exactly $r_n(w)+1$ indices $j$ such that $l_j(w)\geq l_n(w)$.
Since $L(u)$ is a rearrangement of $L(w)$,
there are exactly $r_n(w)+1$ indices $j$ in $u$ such that $l_j(u)\geq l_n(w)$.
Since there are no dots in $D$,
we have $l_i(u)=0$.   Along with the fact that there are no dots in $E$, we see that
 if $j$ is an index
such that $l_j(u)\geq l_n(w)>0$, then the box $(j, u_j)$ lies in either $F$
or $R'$.
Since $r_n(w)+l_n(w)=n-1$, by the assumption $r_n(w)>\frac{n-1}{2}$,
we see that $|\{\text{dots in $R'$}\}|=r_n(w)>l_n(w)$, and so,
for each $1\leq k\leq l_n(w)$, the box $(k, u_k)$ is in $R'$.
Clearly, for $1\leq k\leq l_n(w)$,
we have   $l_k(u)< l_n(w)$.
Hence the maximal possible number of indices $j$ such that $l_j(u)\geq l_n(w)$ is
\[|\{\text{dots in $F$}\}|+|\{\text{dots in $R'$}\}|-l_n(w),\]
which, together with   \eqref{buch},  becomes
\[|\{\text{dots in $R'$}\}|=r_i(u)=r_n(w),\]
contrary to the fact that
there are exactly $r_n(w)+1$ indices $j$ such that $l_j(u)\geq l_n(w)$.
This verifies  $l_n(w)=0$, and hence there are no dots in $F$.

Since  $l_n(w)=0$, we obtain that $w_n=1$.
As there are no dots in $E$ and $F$, we have $i=n$,
and so $r_n(u)=r_n(w)=n-1$, leading to $u_n=1$.
Hence we have $w_n=u_n=1$.

Now we have  proved that either  $w_n=u_n=n$ or $w_n=u_n=1$.
This allows us to finish the proof of the proposition  by induction on $n$.
To be more specific,
if $w_n=u_n=n$, let $w'\in S_{n-1}$ (respectively, $u'\in S_{n-1}$) be
the permutation obtained from $w$ (respectively, $u$)
by ignoring $w_n$ (respectively, $u_n$), while if
$w_n=u_n=1$, let $w'\in S_{n-1}$ (respectively, $u'\in S_{n-1}$) be
the permutation obtained  from $w$ (respectively, $u$)
by ignoring $w_n$ (respectively, $u_n$) and then decreasing each of the
remaining elements by $1$. It is easily seen that   $w'$ and $u'$ satisfy the assumptions
in the proposition.  So, by induction, we have $w'=u'$, which yields $w=u$.
This completes the proof.
\qed

\subsection{Proof of Theorem \ref{main1}}

To present a proof of Theorem \ref{main1}, we still need several lemmas.
Recall that $d(w)$ denotes the inversion code of a permutation $w$.
View
\[d\colon w\longmapsto d(w)\]
as a bijection from $S_n$ to the set $\{(v_1,\ldots,v_n)\colon 0\leq v_i\leq n-i\}$.
With this notation, for any vector $v=(v_1,\ldots,v_n)$ with $0\leq v_i\leq n-i$,
$d^{-1}(v)$ is the corresponding permutation in $S_n$. For $w\in S_n$,
write $w^c=w^c_1w^c_2\cdots w^c_n$ for
the complement of $w$, that is,  $w^c_i=n+1-w_i$ for $1\le i\le n$.
It is routine to check that
\[d(w^c)=\delta_n- d(w),\]
 and so we have
\begin{equation}\label{BU}
w^c=d^{-1}(\delta_n-d(w)).
\end{equation}

\begin{lem}\label{LL-2}
Let $w\in S_n$ be a permutation. Assume that
$x^{\delta_n} \S_w(x^{-1})$ is equal to a Schubert polynomial $\S_{w^\ast}(x)$.
Then
\begin{equation}\label{PP-1}
w^\ast=(d^{-1}(d^{\,t}(w^{-1})))^c,
\end{equation}
which is a permutation in $S_n$.
\end{lem}

\pf By Proposition  \ref{small},
the smallest term
of $\S_w(x)$ in the reverse lexicographic order  is $x^{d^{\,t}(w^{-1})}$.
Hence the leading term of $x^{\delta_n} \S_w(x^{-1})$ is  $x^{\delta_n-d^{\,t}(w^{-1})}$.
On the other hand, by Proposition  \ref{large}, the  leading
 term of $\S_{w^\ast}(x)$ is $x^{d(w^\ast)}$.
 Since $x^{\delta_n} \S_w(x^{-1})=\S_{w^\ast}(x)$, we have
\begin{equation}\label{XY-1}
d(w^\ast)=\delta_n-d^{\,t}(w^{-1}).
\end{equation}
Write $d(w^{-1})=(a_1,a_2,\ldots,a_n)$ and $d^{\,t}(w^{-1})=(b_1,b_2,\ldots,b_n)$.
For  $1\leq i\leq n$, since $0\leq a_i\leq n-i$, we have $0\leq b_i  \leq n-i$.
This, together with  \eqref{XY-1}, implies that $w^\ast$ is a permutation in $S_n$.
Moreover, by \eqref{BU} and  \eqref{XY-1}, we are led to
\[w^\ast=d^{-1}(\delta_n-d^{\,t}(w^{-1}))=(d^{-1}(d^{\,t}(w^{-1})))^c,\]
as desired.
\qed

The second lemma is a well-known  characterization of 132-avoiding permutations,
see for example  \cite[Chapter 1]{Sta}.

\begin{lem}\label{LL-6}
A permutation  $w\in S_n$ is a 132-avoiding permutation  if and
only if its inversion code $d(w)$ is a partition. Moreover, the inverse $w^{-1}$
of $w$ is also
a 132-avoiding permutation with  $d(w^{-1})=d^{\,t}(w)$.
\end{lem}

Based on Proposition \ref{lemma-a} and
Lemma \ref{LL-6}, we are ready to prove the final lemma.

\begin{lem}\label{LL-4}
Let $w$ be a permutation in $S_n$, and let $w^*$
be  as defined in   \eqref{PP-1}.
Then $(w^\ast)^\ast=w$ if and only if $w$ avoids the patterns 132 and 312.
\end{lem}

\pf We first show that if   $w$ avoids  132 and 312, then $(w^\ast)^\ast=w$.
Since $w$ is 132-avoiding, by Lemma \ref{LL-6}, $w^{-1}$ is also 132-avoiding.
By  \eqref{XY-1} and Lemma \ref{LL-6}, we have
\[d(w^\ast)=\delta_n-d^{\,t}(w^{-1})=\delta_n-d(w),\]
yielding  that  $w^\ast=w^c$ by \eqref{BU}.
As $w$ is 312-avoiding, $w^c$ is 132-avoiding, and hence  $w^*$ is 132-avoiding.
Thus  we obtain that
\[(w^\ast)^\ast=(w^c)^\ast=(w^c)^c=w,\]
where the second equality follows by applying
\eqref{XY-1} and Lemma \ref{LL-6}, this
time to the 132-avoiding permutation $w^c$.

Let us proceed to prove the reverse direction. By definition,
\[w=(w^\ast)^\ast=(d^{-1}(d^{\,t}((w^\ast)^{-1})))^c.\]
So we have
$w^c=d^{-1}(d^{\,t}((w^\ast)^{-1})),$
and thus
\begin{equation}\label{disp}
d(w^c) =d^{\,t}((w^\ast)^{-1}).
\end{equation}
 Notice that for any vector $v$ of nonnegative integers,
its transpose $v^t$ as defined in \eqref{trans} is a partition.
Thus, by \eqref{disp} we see that $d(w^c)$ is a partition, and so
it follows from Lemma \ref{LL-6} that  $w^c$
is 132-avoiding, or equivalently, $w$ is 312-avoiding.

We still need to show that $w$ is 132-avoiding. Since $d^{\,t}(w^{-1})$ is a partition,
by Lemma   \ref{LL-6}, there exists a 132-avoiding permutation $u\in S_n$ such that
\begin{equation}\label{PQ-2}
d(u)=d^{\,t}(w^{-1}).
\end{equation}
So we have
\begin{equation}\label{PQ-1}
w^\ast=(d^{-1}(d^{\,t}(w^{-1})))^c=(d^{-1}(d(u)))^c=u^c.
\end{equation}
We now consider $(u^c)^\ast$.
Since  $d^{\,t}((u^c)^{-1})$ is a partition, there exists a 132-avoiding permutation
$v\in S_n$ such that
\begin{equation}\label{PQ-3}
d(v)=d^{\,t}((u^c)^{-1}),
\end{equation}
and so
\[(u^c)^\ast=(d^{-1}(d^{\,t}((u^c)^{-1})))^c=(d^{-1}(d(v)))^c=v^c,\]
 which along with \eqref{PQ-1}
gives $(w^\ast)^\ast=v^c$. By the assumption $(w^\ast)^\ast=w$, we obtain that
\begin{equation}\label{PQ-4}
w=v^c.
\end{equation}

Using \eqref{PQ-2}, \eqref{PQ-3} and \eqref{PQ-4}, we can  prove  the following two claims.

\noindent
Claim 1: $R(w)$ is a rearrangement of $R(u)$. Recall that $R(w)=(r_1(w),\ldots,r_n(w))$, where $r_i(w)$ is equal to the
number of dots lying to the  upper right of the box $(i, w_i)$.
Notice  that in
 the inversion code $d(w)=(d_1,\ldots,d_n)$, the entry $d_i$ equals the number of dots lying  to the  lower left of the box $(i, w_i)$. Moreover,   the matrix representation of $w^{-1}$ is the transpose of the matrix representation of $w$.
Therefore,
$R(w)$ is   a rearrangement of $d(w^{-1})$. For the same reason,
$R(u)$ is  a rearrangement of $d(u^{-1})$.

On the other hand, since $u$ is 132-avoiding, it follows from Lemma \ref{LL-6}
that $d(u)=d^{\,t}(u^{-1})$. Combined with
\eqref{PQ-2}, we have
\[d^{\,t}(u^{-1})=d^{\,t}(w^{-1}).\]
It is easy to check  that for any two vectors $\alpha, \beta\in \mathbb{Z}_{\geq 0}^n$,
if $\alpha^t=\beta^t$, then $\alpha$ is a rearrangement of  $\beta$.
So   $d(w^{-1})$ is a rearrangement of $d(u^{-1})$.
We have explained  that  $R(w)$ is   a rearrangement of $d(w^{-1})$ and
$R(u)$ is  a rearrangement of $d(u^{-1})$. Hence
 $R(w)$ is a rearrangement of $R(u)$. This proves Claim 1.

\noindent
Claim 2: $L(w)$ is a rearrangement of $L(u)$. Since $v$ is 132-avoiding, by Lemma \ref{LL-6}
we have $d(v)=d^{\,t}(v^{-1})$. Moreover, by \eqref{PQ-4} we see that $v^{-1}=(w^c)^{-1}$.
Hence,
\[d(v)=d^{\,t}(v^{-1})=d^{\,t}((w^c)^{-1}).\]
In view of \eqref{PQ-3}, we get
\[d^{\,t}((w^c)^{-1})=d^{\,t}((u^c)^{-1}).\]
By the same arguments as in Claim 1,   $R(w^c)$ is a rearrangement of $R(u^c)$.
Noticing that $L(w)=R(w^c)$ and $L(u)=R(u^c)$, we conclude that
$L(w)$ is a rearrangement of $L(u)$. This verifies Claim 2.

Since $w$ is 312-avoiding and $u$ is 132-avoiding, combining
 Proposition  \ref{lemma-a}, Claim 1 and Claim 2, we obtain
that $w=u$, and so $w$ is  132-avoiding. This finishes  the proof.
\qed

We are now in the position  to provide  a proof of Theorem \ref{main1}.

\noindent
{\it Proof of Theorem \ref{main1}.}
We first prove the sufficiency.
Assume that $w\in S_n$ avoids the patterns 132 and 312. Since $w$ is 132-avoiding,
it follows from \cite[Chapter IV]{Mac} that $\S_w(x)=x^{d(w)}$.
So we have \begin{equation}\label{Final-1}
x^{\delta_n} \S_w(x^{-1})=x^{\delta_n-d(w)}.
\end{equation}
As $\delta_n-d(w)$
is the inversion code of $w^c$, we see that $x^{\delta_n-d(w)}=x^{d(w^c)}$.
Since $w$ is 312-avoiding, $w^c$ is 132-avoiding, and so
 we have
 \begin{equation*}
\S_{w^c}(x)=x^{d(w^c)}=x^{\delta_n-d(w)}.
\end{equation*}
This together with \eqref{Final-1} yields that $x^{\delta_n} \S_w(x^{-1})$ is the Schubert
polynomial $\S_{w^c}(x)$.

It remains to prove the necessity. Assume that $x^{\delta_n} \S_w(x^{-1})$ equals the Schubert polynomial
$\S_{w^\ast}(x)$. Equivalently, we have  $\S_w(x)=x^{\delta_n} \S_{w^\ast}(x^{-1})$,  which
along with
Lemma \ref{LL-2} leads to
$w=(w^\ast)^{\ast}$. Invoking  Lemma \ref{LL-4},
we conclude that $w$ avoids  132 and $312$. This completes the proof.
\qed

\section{Concluding remarks on Conjecture \ref{cytj}}\label{TEE}

The necessity of Conjecture \ref{cytj} predicts that if $x^{\mu}\S_w(x^{-1})$
is  a Schubert polynomial,
then $w$ avoids the  patterns  1432, 13254, 14253, 24153, 31524, 361452.
In the case when $w$ is the pattern 1432, we  explain that  $x^{\mu}\S_{1432}(x^{-1})$
cannot be a Schubert polynomial.
Notice that
\begin{equation}\label{PP9}
x^{\mu}\S_{1432}(x^{-1})=\frac{x^{\mu}}{x_1^2x_2^2x_3}(x_1^2+x_1x_2+x_1x_3+x_2^2+x_2x_3).
\end{equation}
Suppose otherwise that
$x^{\mu}\S_{1432}(x^{-1})$ is a Schubert polynomial, say,
$x^{\mu}\S_{1432}(x^{-1})=\S_{w'}(x)$ for some permutation $w'\in S_n$.
In the reverse lexicographic order,  it follows from \eqref{PP9} that
$\frac{x^{\mu}}{x_1^2x_2^2x_3}\,x_2x_3$ is the leading monomial of $\S_{w'}(x)$.
Let $d(w')=(d_1,\ldots,d_n)$ be the inversion code of $w'$. By Proposition \ref{large},
\begin{equation}\label{PP8}
d_1=\mu_1-2,\ d_2=\mu_2-1,\ d_3=\mu_3,\ d_4=\mu_4,\ldots, d_n=\mu_n.
\end{equation}
By Theorem \ref{mt1}, $\mu$ is a partition. So we have
\begin{align}\label{bd}
d_1+2\ge d_2+1\ge d_3\ge d_4\ge\cdots\geq d_n.
\end{align}

On the other hand, we see from  \eqref{PP9} that
$\frac{x^{\mu}}{x_1^2x_2^2x_3}\,x_1^2$
 is the smallest monomial of $\S_w(x)$ in the reverse lexicographic order.
 By \eqref{PP8},
 \begin{equation}\label{PP7}
\frac{x^{\mu}}{x_1^2x_2^2x_3}\,x_1^2=x_1^{d_1+2}x_2^{d_2-1}x_3^{d_3-1}x_4^{d_4}\cdots x_n^{d_n}.
\end{equation}
 Denote by $D_{\mathrm{bot}}$  the bottom RC-graph of $w'$, which  contains the first $d_i$ boxes in row $i$.
From Section \ref{ma1}, we know that the monomial \eqref{PP7} can be generated from $D_{\mathrm{bot}}$ by applying
 ladder moves. To be specific, we need to move two crosses from
 $D_{\mathrm{bot}}$ to the first row:   one  from  the second row and   one
 from the third row. This, in view of \eqref{bd}, ensures that $d_2=d_1+1$ and $d_3=d_2+1$.
 However, in this case we can apply ladders moves to generate a monomial
 \[x_1^{d_1+2}x_2^{d_2}x_3^{d_3-2}x_4^{d_4}\cdots x_n^{d_n},\]
 which is smaller than the monomial in  \eqref{PP7}, see Figure \ref{P100} for an illustration.
  \begin{figure}[h]
\setlength{\unitlength}{1.5mm}
\begin{center}
\begin{picture}(200,6)

\ytableausetup{boxsize=.95em}
\ytableausetup{aligntableaux=center}

\begin{ytableau}
\scriptstyle$+$ &  \scriptstyle & \scriptstyle &\scriptstyle &\scriptstyle&\scriptstyle \\
\scriptstyle$+$ &  \scriptstyle$+$ & \scriptstyle &\scriptstyle&\scriptstyle \\
\scriptstyle$+$ &  \scriptstyle$+$ &\scriptstyle $+$ &\scriptstyle
\end{ytableau}
 $\rightarrow$

\begin{ytableau}
\scriptstyle$+$ &  \scriptstyle & \scriptstyle$+$ &\scriptstyle &\scriptstyle&\scriptstyle \\
\scriptstyle$+$ &  \scriptstyle & \scriptstyle &\scriptstyle&\scriptstyle \\
\scriptstyle$+$ &  \scriptstyle$+$ &\scriptstyle $+$ &\scriptstyle
\end{ytableau}
 $\rightarrow$

\begin{ytableau}
\scriptstyle$+$ &  \scriptstyle & \scriptstyle$+$ &\scriptstyle &\scriptstyle&\scriptstyle \\
\scriptstyle$+$ &  \scriptstyle & \scriptstyle &\scriptstyle$+$&\scriptstyle \\
\scriptstyle$+$ &  \scriptstyle$+$ &\scriptstyle  &\scriptstyle
\end{ytableau}
 $\rightarrow$

\begin{ytableau}
\scriptstyle$+$ &  \scriptstyle & \scriptstyle$+$ &\scriptstyle &\scriptstyle$+$&\scriptstyle \\
\scriptstyle$+$ &  \scriptstyle & \scriptstyle &\scriptstyle&\scriptstyle \\
\scriptstyle$+$ &  \scriptstyle$+$ &\scriptstyle  &\scriptstyle
\end{ytableau}
 $\rightarrow$

\begin{ytableau}
\scriptstyle$+$ &  \scriptstyle & \scriptstyle$+$ &\scriptstyle &\scriptstyle$+$&\scriptstyle \\
\scriptstyle$+$ &  \scriptstyle & \scriptstyle $+$ &\scriptstyle&\scriptstyle \\
\scriptstyle$+$ &  \scriptstyle &\scriptstyle  &\scriptstyle
\end{ytableau}

\end{picture}
\end{center}
\vspace{.2cm}
\caption{ Ladder moves in the case when $d_2=d_1+1$ and $d_3=d_2+1$.}
\label{P100}
\end{figure}
Hence, the assumption that
$x^{\mu}\S_{1432}(x^{-1})$ is a Schubert polynomial is false.

From  the above arguments, we see that even in the special case of $w=1432$,
the proof is not trivial.
We would like to mention that Fink,  M\'esz\'aros and  St.\,Dizier
\cite{Fin}  showed that if $\sigma$ is a pattern of $w$, then
$\S_w(x)$ can be expressed as a monomial times $\S_\sigma(x)$
(in reindexed variables) plus a polynomial with nonnegative coefficients.
It seems  that this result might  be helpful in the
investigation  of Conjecture  \ref{cytj}.

\footnotesize{

N{\scriptsize EIL} J.Y. F{\scriptsize AN}, D{\scriptsize EPARTMENT OF} M{\scriptsize ATHEMATICS}, S{\scriptsize ICHUAN} U{\scriptsize NIVERSITY}, C{\scriptsize HENGDU} 610064, P.R. C{\scriptsize HINA.} Email address: fan@scu.edu.cn

P{\scriptsize ETER} L. G{\scriptsize UO}, C{\scriptsize ENTER FOR} C{\scriptsize OMBINATORICS}, N{\scriptsize ANKAI} U{\scriptsize NIVERSITY}, T{\scriptsize IANJIN} 300071, P.R. C{\scriptsize HINA.}  Email address: lguo@nankai.edu.cn

N{\scriptsize ICOLAS} Y. L{\scriptsize IU}, C{\scriptsize ENTER FOR} C{\scriptsize OMBINATORICS}, N{\scriptsize ANKAI} U{\scriptsize NIVERSITY}, T{\scriptsize IANJIN} 300071, P.R. C{\scriptsize HINA.}  Email address: yiliu@mail.nankai.edu.cn

}

\end{document}